# Nonorthogonal Bases and Phase Decomposition: Properties and Applications.[1]

**Sossio Vergara**
ITI Cardano, Rome, Italy

## Abstract

In a previous paper (Vergara, S., *On Generic Frequency Decomposition. Part 1: Vectorial Decomposition*. Dig. Sig. Proc. Vol. 17 N. 2 2007) it was discussed the viability of functional analysis using as a basis a couple of generic functions, and hence vectorial decomposition. Here we complete the paradigm exploiting one of the analysis methodologies developed there, but applied to phase coordinates, so needing only one function as a basis. It will be shown that, thanks to the novel iterative analysis, any function satisfying a rather loose requisite is ontologically a basis. This in turn generalizes the polar version of the Fourier theorem to an ample class of nonorthogonal bases. The main advantage of this generalization is that it inherits some of the properties of the original Fourier theorem. As a result the new transform has a wide range of applications and some remarkable consequences. The new tool will be compared with wavelets and frames. Examples of analysis and reconstruction of functions using the developed algorithms and generic bases will be given. Some of the properties, and applications that can promptly benefit from the theory, will be discussed. The implementation of a matched filter for noise suppression will be used as an example of the potential of the theory.

Keywords: Fourier analysis, wavelets, nonorthogonal, polar.

vsossio@gmail.com

---

## 1. Introduction.

The Fourier theorem is one of the most valuable mathematical tool, it is used in all sort of applications, from mobile phones to system theory. In its standard form the theorem states the possibility to reconstruct a function using a series of sines and cosines. Over the years different generalizations have been devised for the theorem. One of the most useful is the use of bases other than the orthogonal pair, for example using wavelets and frames. Alas these tools for nonorthogonal bases sport a high computational complexity, tolerate only few specially built functions as bases, and require different bases for analysis and reconstruction (biorthogonal and dual bases) [6, 7, 8, 9]. These constraints have limited the diffusion of the wavelets and frames to special applications, compared to the pervasiveness of the original Fourier theorem.

A recent development is the Common Waveform Analysis [10, 14, 15]. There, a couple of even and odd special functions as the square wave, triangular and the like constitutes the basis. However, due to the limitation of the mathematics involved (that is based on an inner product computation), this theory still requires the use of biorthogonal bases, works only with a limited number of special bases and vectorial decomposition.

A different approach has been introduced in [1]. There, two procedures were developed for the analysis. One method has been called "direct" or brute force and requires the solution of a system of equations (much like the algorithm used in frame theory). The second method has been dubbed "indirect" because uses a novel iterative algorithm. The two methods have been used together in [1] to validate each other, although it has been briefly illustrated the superiority of the indirect method over the direct one.

The theory developed in [1], compared to the known tools as the wavelets, the frames and the common waveform analysis, admits a much larger class of functions as bases and, most notably, utilizes the same basis for analysis and reconstruction. All thanks just to a change in the analysis paradigm. The idea is that there is no reason why we should limit ourselves to the use of the inner product in the analysis. In effect the definition of basis does not mention the analysis phase, requiring only the possibility of reconstruction of any function of the given space in terms of a combination of the basis (plus the usual constraints of unicity and convergence of the reconstruction). In other words we are left free to choose the analysis method (here we prefer the word "decomposition" for reasons that will be clear below). Traditionally this freedom is not exploited as the vast majority of the established algorithms for nonorthogonal bases employ the same inner product computation of the original Fourier theorem. And thus a "direct" approach. The advantage is the possibility of calculating any of the components independently, but at the same time, oblige to use biorthogonal functions when the basis is nonorthogonal and greatly limits the choice of the basis. Instead the recursive computation methodology (the "indirect" method) exploits an iterative change of coordinates between the orthogonal basis and the and the new, nonorthogonal, one [1].

To find the components of a signal with the indirect methodology, one has to start with the decomposition of the signal in the usual orthogonal basis. This generates a representation in the Hilbert space. Then an iterative algorithm for change of reference translate these components from the Hilbert space to the new basis. The originality of the method is that we do not try to find directly the components of a function on a nonorthogonal basis. Instead, we switch from a known representation (the Fourier components in the Hilbert space) to a new, equivalent, one. For the transitive property of equality, if one is a representation of the signal so is the other.

In literature there is another main iterative algorithm for function analysis: the Empirical Mode Decomposition (EMD), as used for example in Hilbert-Huang transform [17]. However there are substantial differences between this approach and the EMD. For one EMD is, by definition, empirical: the basis functions are derived from the data. Instead here, one first chooses the basis, then the parameters of the decomposition are computed. The benefit is that one can pick the basis that is best suited to a particular problem.

Matching Pursuit (MP) [12] is another common algorithm for function analysis. There a dictionary, generally consisting of a large collection of time-frequency atoms, is employed in the search for the best sparse representation of a signal, minimizing the error. The tool is most useful for compression and coding but is very computing intensive. Instead in [1] there is much more freedom in the choice of the basis, and the decomposition is purely in frequency. And for each basis an exact representation of the signal is computed up to any chosen frequency. In other words the error can be confined at highest frequencies and with a very efficient algorithm.

An added advantage of this method is that one can still use the metric of the Hilbert space to assert the convergence of the procedure. The only drawback of the iterative method is that to find the component at a given frequency, generally all the components at lower frequency must be computed. But we think that this is not an issue given the possibilities that the new theory opens up. A further advantage is that one can now use the same basis for analysis and reconstruction (in contrast with the usual tools for nonorthogonal bases).

But probably the main benefit of the "indirect methodology" is that it is applicable even to polar decomposition and this will help us to surmount the borders of linear dependence.

As a matter of fact all of the cited tools are based on the methods of linear algebra, so inherently vectorial. However, at the time of interpreting the analysis results, these tools are no match for the simplicity of a Fourier power spectrum in terms of amplitude and phase.

The natural extrapolation of the theory developed in [1] is thus applying the same computing method to phase decomposition. The benefit of phase (in orthogonal terms: "polar") decomposition is in dealing with a single function as a basis. And a single function as a basis has the indisputable advantage (in contrast with the vectorial tools) that some of the properties of the Fourier theorem can be extended also to nonorthogonal bases, and this greatly enhances the applicability of the tool, as it will be clear in the following.

When orthogonal bases are involved, polar and vectorial representations are essentially the same thing, as there is a trivial equation connecting the two coordinate systems. Instead, with nonorthogonal bases, the vectorial and phase decompositions of a function are completely different beasts, and there is no simple way to pass from one to the other. An example could clarify the point. Imagine of having as a function a square wave with arbitrary phase. If we use vectorial decomposition and a basis consisting of even and odd square waves, a viable basis according to [1], [10], [14], [15] (they all give the same results, although the last papers are based on traditional approach and biorthogonal bases), it is evident that we would need an infinite series of square waves to reconstruct the function. Because the nonorthogonal vectorial decomposition cannot easily characterize the arbitrary phase. Instead, when using phase decomposition and a single square wave as a basis [13], the result of the analysis is a single couple of parameters at the given frequency: amplitude and phase. Much more efficient and comprehensible.

The disadvantage of the phase decomposition is that its outcome is a set of two parameters: amplitude and phase (or more generally shift) that are not homogeneous, differing dimensionally, and hence preventing the use of matrices and linear algebra in the computations (that would be precluded anyway because in case of nonorthogonal bases the resulting systems will be nonlinear, as it will be shown below). As a consequence, for orthogonal bases the vectorial analysis is preferred, as in the most common flavors of the Fourier Theorem. Whereas, when using nonorthogonal bases, the phase decomposition is more widely applicable and delivers more interesting results, even with the added burden of dealing with couples of non-homogeneous parameters.

A previous paper [13] introduced the phase decomposition over nonorthogonal bases but with a focus on a special application. There, it was demonstrated that the square wave is one of the viable bases for phase decomposition. As the square wave is the natural output of digital systems, it was exploited in the design of very efficient, multiplierless, signal synthesizers. The systems employing the square wave are very frugal on computing demands and suitable for many applications.

The goal of this article is to disclose the bigger picture, revealing some of the properties and consequences of the theory of phase decomposition over nonorthogonal bases and indicating other applications.

## *2. The iterative analysis methodology.*

Here the rationale behind the new computation scheme will be briefly summarized in order to introduce a fast analysis algorithm and other consequences of the decomposition.

Lemma:
Given an Hilbert space **H** with the usual orthonormal basis, then any function $S(x) \in \mathbf{H}$ spans the space when using the same frequency-phase reconstruction algorithm of the polar Fourier decomposition. I.e. $\{S(nx)\}$ is complete (*n* being the frequency).

We shall prove the feasibility for real periodic functions $f(x)$, $S(x) \in L^2[-\pi, +\pi]$ (the space of periodic Lebesgue square integrable functions). While the extension to complex valued functions, different periods and transforms is straightforward.

Given any periodic function $f(x) \in L^2[-\pi, +\pi]$ satisfying the Dirichlet conditions, it can be expressed as a Fourier series. We omit here an eventual average (a DC component) from the series as it is a simple constant that will not change our conclusions:

$$f(x) = \sum_{k=1}^{\infty} b_k \cos(kx + \vartheta_k) \tag{1}$$

and given another nonzero periodic function with Fourier series:

$$S(x) = \sum_{p=1}^{\infty} s_p \cos(px + \phi_p) \tag{2}$$

It is easy to verify that one can shift and amplitude scale the $S(nx)$ (at any frequency $n \in N$) in order match any cosine component, using two parameters $M$, $\Theta$ for amplitude and phase:

$$a \cdot \cos(nx + \vartheta) + error = M \cdot s_1 \cos(nx + \phi_1 + \Theta) + M \sum_{p=2}^{\infty} s_p \cos\left[n\left(p(x+\Theta) + \phi_p\right)\right] \tag{3}$$

The last term is the error at higher frequencies due to the nonorthogonality of the $S(x)$. In the procedure we use the error function, i.e. the $f(x)$ minus the previous reconstruction. At n = 1 the error function $f_{e1}(x)$ is $f(x)$ itself. At n=2 we require that the harmonic at frequency one is cancelled from the error function (the term harmonic here refers only to a Fourier, sinusoidal, component of a function). The resulting equation being:

$$f_{e2}(x) = f_{e1}(x) - M_1 S(x + \Theta_1) = \sum_k b_k \cos(kx + \vartheta_k) - M_1 \sum_p s_p \cos(p(x+\Theta_1) + \phi_p) \tag{4}$$

One can compute the unknown parameters from the equations:

$$M_1 = \frac{b_1}{s_1} \tag{5}$$

$$\Theta_1 = \vartheta_1 - \phi_1 \tag{6}$$

The $f_{e2}(x)$ with the computed parameters has only harmonics from frequency two up. One can then use the $f_{e2}(x)$ as the new function iterating the (4) at frequency two with $M_2 S(2x+\Theta_2)$ to find $M_2$ and $\Theta_2$ using again the corresponding equations (5) and (6). The final result being the recursive equation:

$$f_{eN+1}(x) = f_{eN}(x) - M_N S(Nx + \Theta_N) = f(x) - \sum_{i=1}^{N} M_i S(ix + \Theta_i) \tag{7}$$

In (7) the lowest harmonic (at any frequency N) vanishes from the error function for each cycle of the iteration. The iteration can be stopped after a given number of cycles (components) or when the norm of error function (7) results lower than a threshold.
I.e. any $\{S(nx)\} \in H$ with $n=1,2,3,…, \infty$, spans the Hilbert space when using this iterative procedure to determine the coefficients.

Please note that we only need to define the antitransform, or synthesis, algorithm (the rightmost term in the (7) ). The transform (or analysis) is derived from it.
The entire process becomes somehow the "*de-construction*" of the function in the given basis.
The (5), (6) and (7) form an iterative reference change from the sinusoidal basis (in the Hilbert space) to the new basis *S(nx)*. The algorithm works assembling a reconstruction in terms of *S(nx)* having the same Fourier spectrum as the original function. And, thanks to the unicity of the Fourier theorem, if one represents the *f(x)*, so does the other.
Note that up to now the only requirement on *S(x)* is to have $s_1 \neq 0$ or equivalently a frequency not higher than that of the *f(x)* and integral zero over a period. In order to check whether *S(x)* is a true basis, we need to verify under which conditions the series converges to *f(x)*. Of the different definitions of convergence: pointwise, uniform and norm, the latter appears the most useful. So we can express the previous statement more rigorously requiring that the norm of the error function tends to zero as $N \to \infty$ (following the definition of metrics in Hilbert spaces), or:

$$\lim_{N \to \infty} \| f_{eN}(x) \| = \lim_{N \to \infty} \| f(x) - \sum_{n=1}^{N} M_n S(nx + \Theta_n) \| = 0 \qquad (8)$$

It is easy to prove that *S(x)* is ontologically a basis of the $L^2$ space if most of its energy is concentrated at the fundamental harmonic (the proof is omitted here but can be found in [13]), i.e.:

$$s_1^2 > \sum_{n=2}^{\infty} (s_n)^2 \qquad \text{Convergence Requirement} \qquad (9)$$

In case the right hand side of the inequality is zero, the basis is a pure sinusoid, and we fall back into the original Fourier theorem; as expected from a true generalization.
In the above lemma we used the same waveform *S(x)* at each frequency for simplicity, but it is evident that the proof is equally valid for different functions, one for each frequency (i.e. any $\{S_n(nx)\}$, $n = 1, 2, 3,…, \infty$ is complete). This way even functions that do not satisfy the convergence requirement (9) can be used in the lower frequency part of a basis. For the convergence of the (8) it is sufficient that at higher frequencies we employ a real basis (a function satisfying the (9) ).
We emphasize here that these results, according to the terminology introduced in [12], can also be viewed as a special case of atomic decomposition over a complete dictionary, but can be easily extended to overcomplete dictionaries.

## *3. A Fast Analysis Algorithm.*

The above procedure makes use of the *f(x)* in the computation. However a fast analysis algorithm has been developed from the (7), working directly in Fourier (sinusoidal) frequency domain.
We start from a nonorthogonal basis *S(nx)* having Fourier components as in (3). We then define an amplitude scaled and phase shifted version *S'(nx)* of the original function *S(nx)* such that it has a simple cosine with phase zero and amplitude one as fundamental:

$$S'(nx) = \frac{1}{s_1} S(n(x - \phi_1)) = \sum_{k=1}^{\infty} \frac{s_k}{s_1} \cos(n(k(x - \phi_1) + \phi_k)) \qquad (10)$$

It is always possible doing so and a more compact equation will result at the end.
Switching to the continuous we can write it as *S'(ωx)*, *ω* being a real representing the frequency.
Then *S'(ωx)* has a normalized cosine with phase zero as fundamental at frequency *ω*.
If a set of $\{S'(\omega x)\}$ ($0 < \omega < \infty$) can be used to reconstruct any collection of harmonics (see above) then it shall exist a set $\{M(\omega), \Phi(\omega)\}$ such that we can define an antitransform $T_{S'}^{-1}$ as:

$$f(x) = T_{S'}^{-1}\left[ M(\omega), \Phi(\omega) \right]_{\omega>0}^{\omega\to\infty} \equiv \lim_{\omega_1 \to \infty} \int_{\omega>0}^{\omega_1} M(\omega) S'(\omega x + \Phi(\omega)) d\omega \tag{11}$$

The right hand side of the (11) is derived from the usual Fourier antitransform algorithm in polar form, except that we used here the basis $S'(\omega x)$ (cosine plus higher frequency noise) instead of the pure cosine. Equivalently we can write:

$$\lim_{\omega_1 \to \infty} f(x) - T_{S'}^{-1}[M(\omega), \Phi(\omega)]_{\omega>0}^{\omega=\omega_1} = \lim_{\omega_1 \to \infty} f_{e(\omega<\omega_1)}(x) = 0 \tag{12}$$

Starting from the antitransform (11), we can define the new generalized transform in recursive form, as the Fourier transform, computed at the frequency $\omega_1$, of the difference of $f(x)$ and the reconstruction (12) up to $\omega_1$ (i.e. the Fourier Transform of the error function at frequency $\omega_1$):

$$T_{S'}[f(x)]_{\omega_1} = T_{S'}[f_{e(\omega<\omega_1)}(x)]_{\omega_1} = T_F\left(f(x) - T_{S'}^{-1}[M(\omega), \Phi(\omega)]_{\omega<\omega_1}\right)_{\omega_1} = [M(\omega_1), \Phi(\omega_1)] \tag{13}$$

Notice how, employing the recursive algorithm (13), we reduced the computation of the components in the new basis $S(x)$ to the manipulation of plain Fourier Transforms. Only the definition of the antitransform (11) is needed.
Here $T_{S'}[f(x)]_{\omega_1}$ is the component of $f(x)$ in the basis $S'(\omega x)$ computed at the frequency $\omega_1$, $f_{e(\omega<\omega_1)}$ is the error function at frequency $\omega_1$.
$T_F(\ldots)_{\omega_1}$ is the Fourier transform in polar coordinates computed at $\omega_1$. We can here employ explicitly amplitude and phase of the Fourier transform because we switched from $S(x)$ to the $S'(x)$.
At the lowest frequency $\omega_0$ the error function is the function itself (there is no reconstruction yet), so that the $M(\omega_0)$, $\Phi(\omega_0)$ are the very Fourier parameters of $f(x)$. Then we step up the frequency and recompute the (13). The (13) includes uniquely Fourier transforms and thus ensures that any function that has a Fourier reconstruction, also has a reconstruction in terms of the new basis $S'(x)$. The (13) can be computed directly in the Fourier frequency-phase domain.
Once computed the complete set of coefficients $\{M(\omega), \Phi(\omega)\}$ one can scale them to reflect the original $S(x)$. In Fig. 3 there is the pseudocode of a fast analysis algorithm.
Here we used the letter M from the word "modulus" for the first coefficient as a vestige of its orthogonal ancestor, although we shall use also the more correct term "amplitude" for this coefficient, and likewise for the frequency we shall use the letter $\nu$ or the $\omega$ when necessary to remark the derivation from the Fourier theorem. The complexity of the analysis algorithm depends on the number and distribution of the harmonics of the basis $S(x)$ but from the flowchart in Fig. 3 it results that in the worst case, for an N point analysis, the complexity is at most:

$$O(N \log(N)) + O\left( \sum_{i=1}^{N} \frac{N}{i} \right) \approx 2 \times O(N \log(N)) \tag{14}$$

One first needs the FFT of the basis but it is done just once. Then one needs to compute the FFT of the function (the first term above) and then O(N) operations to calculate the first component of the basis and its effect on the error function. To compute the component of the basis at frequency two and its effect on the error function one needs O(N/2) operations and so forth for the remaining components up to frequency N. This is the meaning of the second term in the (14). So at the end just about the double of the FFT.

## 4. Results and Properties of the New Decomposition.

To illustrate the capabilities of the new transform in Fig.1 one period of a function is plotted together with its reconstruction, using a generic basis (nonorthogonal and asymmetric), and limiting the reconstruction to ten components. The original function is dotted while the re-synthesized function is solid, in the lower left part of the figure there is a reduced plot of the basis. In Fig. 2 the same function with the same basis has been plotted, but the approximation is extended up to thirty components. The noise is now lower and confined to higher frequency as reasonably expected from a converging basis. The meaning of these figures is to show how a periodic signal can be approximated with a series of another nonorthogonal and asymmetric function (thus being a real basis).

If we instead analyze and reconstruct a function using a "*diverging basis*" (we use this oxymoron for a function for which the (9) does not hold), the noise will increase with the order N of the series: the function $S(x)$ still spans the space but is not a well-behaved basis.

The Fourier Transform has some nice properties, some of these hold true also for the new one, while others need to be changed to take into account the nonorthogonality of the basis. For lack of space we can only introduce here some of the main results just to reveal the overall framework of the theory.

A property that is applicable also to the new Transform is the proportionality [2] [4]. For the Fourier Transform if:

$$F(\omega) = T_F(f(x)) \qquad (15)$$

Then:

$$T_F(k f(x)) = k T_F(f(x)) = k F(\omega) \qquad (16)$$

And the same is true for a generic basis. But first we need to rewrite the previous equation in phase coordinates, and the result is a multiplication of the amplitudes while the phases are unchanged. So if:

$$T_S(f(x)) = \{M(\omega), \Phi(\omega)\} \qquad (17)$$

Then it is straightforward to verify from the (13):

$$T_S(k f(x)) = \{k M(\omega), \Phi(\omega)\} \quad , \quad T_S^{-1}\{k M(\omega), \Phi(\omega)\} = k f(x) \qquad (18)$$

The (time) translation also holds, in this case only the phase will be affected (see again the (13)):

$$T_S(f(x - x_0)) = \{M(\omega), \Phi(\omega) - \omega x_0\} \quad , \quad T_S^{-1}\{M(\omega), \Phi(\omega) - \omega x_0\} = f(x - x_0) \qquad (19)$$

These last two equations are exactly the same as in the polar Fourier transform.

Instead it is easy to verify that in general, due to the nonorthogonality of the basis $S$, the linearity does not hold:

$$T_S^{-1}[T_S(f(x) + g(x))] \neq T_S^{-1}[T_S(f(x)) + T_S(g(x))] \qquad (20)$$

To be convinced one can think to the sum of two square waves at the same frequency but with different phases that is generally not a square wave, thus having a complex decomposition in the square wave basis, whereas each original square wave has only one component[2].

Although in certain important cases the (20) can be instead an equation and this fact will be exploited in a test application in the last section. The (20) behaves linearly (is an equation) if the components of $f(x)$ and $g(x)$ in the S-basis frequency domain are disjoint (let us call them S-orthogonal) defined as:

$$f(x) \perp_S g(x) \Leftrightarrow \sum_i M_{fi} M_{gi} = 0 \qquad (21)$$

Where $\{M_{fi}\}$ is the amplitude spectrum of $f(x)$ in the basis $S(x)$ and similarly $\{M_{gi}\}$ is for $g(x)$. Here as in most part of the paper we, for clarity, used the series; however the formulas are easily translated to the transforms and the integrals.

Please note that the porting of these properties of the Fourier theorem to the nonorthogonal generalization is possible only thanks to the choice of phase decomposition. For the vectorial decomposition, as is used in the most common analysis tools, it would be impossible. And are exactly these properties that make the polar generalization much more powerful than any vectorial tool.

---

[2] The sum by definition must be a linear operation in both domains of (20). Instead the sum of two copies of a nonorthogonal basis at the same frequency in time domain, always generates a new waveform (thus having a complex spectrum in the same basis). And it can not be the result of a linear operation in the frequency domain.

A straightforward application of the theory is in building generalized filters. If we analyze *f(x)* using a basis *S(x)*, we get two functions of the "generic frequency": the amplitude M(ν) and the phase Φ(ν); if we then multiply M(ν) with a transfer function G(ν) and then reconstruct the resulting signal via an inverse transform (11), what we get at the end is a filter.

*f(x)* →Analysis via *S(x)*→[**M(ν), Φ(ν)**]→[**G(ν)M(ν), Φ(ν)**]→Reconstruction via *S(x)*→*y(x)*     (22)

The (22) is our usual definition of a filter, except for the generic basis instead of the sinusoid. In general we can manipulate the phase too but we skip this case for the moment. From the properties (18), (19), (20) above, the algorithm (22) models a class of time-invariant nonlinear systems and, most important, with the basis as *eigenfunction.*

The nonorthogonality of the basis *S(nx)* leads to the inapplicability of the convolution theorem; as it is it works only for the orthogonal Fourier basis. Nevertheless we can still define the convolution in the usual way in the specific frequency domain as:

$$f(t)\,conv_S\,g(t) = T_S^{-1}[T_S(f(t)) conv_S T_S(g(t))] \equiv T_S^{-1}\left[\int_{\nu>0}^{\infty} F_S(\nu)G_S(\nu)S[\nu t + \Theta_S(\nu) + \Gamma_S(\nu)]d\nu\right] \quad (23)$$

where {$F_S(\nu), \Theta_S(\nu)$} is the transform of *f(x)* in the basis *S(x)* and likewise {$G_S(\nu), \Gamma_S(\nu)$} is for *g(t)*. The (23) is simply the explicit form of the (22) including the phase manipulation and is the direct transposition of the usual algorithm of the filter in polar coordinates to the new basis *S(x)*.

The porting of the convolution theorem to the new theory could be theoretically possible, but only including explicitly the basis even in the time domain, and that ends-up being a very involute nonlinear equation, strongly depending on the basis. On the other hand, it means that even simple operations in the new, basis specific, frequency domain could result in very complex, nonlinear time domain transformations, and this fact can be exploited in Non Linear Signal Processing [3].

The importance of this result should not be underestimated. As the sinusoid is the *eigenfunction* of linear systems (and the reason why the Fourier transform is employed in the study of those systems), so the new analysis tool generalizes the mathematical modeling also to nonlinear systems. Thus each basis creates its own specific frequency domain as a result of a transform from the time domain. And equivalently, any nonlinearity has a specific *eigenfunction* that can be used as a basis to model the system.

To test the possibilities of the new tool in Fig. 4 a nonorthogonal basis is plotted. This basis has been generated starting from a sinusoid undergoing a nonlinear amplification. The choice is not casual as the plot in Fig. 4 is very similar to a Jacobi elliptic function that is *eigenfunction* of a well-known nonlinear system: the Duffing oscillator [18].

When we analyze a sinusoid using the function in Fig. 4 as a basis, we get a collection of amplitudes and phases: {M(ν), Φ(ν)}. The amplitudes of the first few components are:
M(1)= 1,1183
M(3)= 0,1486
M(5)= 0,0212
M(7)= 0,0070
Etc.

As the basis is quite similar to the sinusoid at the input, the amplitude converges rapidly to zero. In Fig. 5 there is the plot of the output of the related system for different, very crude, transfer functions. Essentially we assemble an ideal low pass filter that cuts off a variable number of components (see Eq. 22). If one uses all the components (an all-pass filter), he/she gets of course the same sinusoid as the input, and that is the solid curve in Fig 5. The other two curves correspond to fewer components in the antitransform. Evidently when one uses the basis as input, he/she gets at the output that same basis (see again Eq. 22), so the basis is also the *eigenfunction* of the system. With this very simple model and sinusoidal input, just acting on the transfer function, one can get at the output different waveforms.

Of course this is a simplified nonlinear system. In real systems things can get much more complicated. Nonetheless the developed tool can be of much help in the understanding of these phenomena. A more detailed discussion has to be postponed to a forthcoming article.

This finally takes us to generalize the faithful Nyquist theorem. As a matter of fact the sampling rate depends on the basis. Let us imagine a square wave in our old Fourier domain, there the sampling rate of the square wave extends to infinite, while for a sinusoid it is just the double its frequency. But this is true only in case the basis is a sinusoid. If we switch to the frequency domain of the square wave basis, see for example [13], the opposite is true: the sampling rate of a sinusoid is infinite, but that of the square wave needs only to be the double of its frequency. Because the representation of a square wave in the square wave frequency domain is a single component.

In conclusion, with the new deconstruction methodology, any basis defines a class of systems, has its own frequency domain, and the manipulation of a function in a specific frequency domain can be an efficient way to model nonlinear systems.

## 5. Applications

The theory introduced here has many mathematical and physical implications, some of them will be examined in future papers. We now focus just on few possible applications concerning DSP. The main point of the new transform is that almost any function can be used as a basis. And a greater freedom in the choice of the basis can open up a whole new class of applications. A first straightforward application of the theory has been illustrated in [13], where an efficient, multiplierless, signal synthesizer based on square waves has been described.

But many other applications can be easily foreseen, for example in the field of data compression. In this case the basis can be tailored in order to achieve the highest compression ratio, exploiting the characteristics of the data.

The topic is easily extended to data transmission and compression of images. The current methods involve the use of DFT, DCT or DWT (Discrete Transform of Fourier, Cosine or Wavelet types). As a matter of fact the most promising approach to image compression is the use of Haar wavelets [5][16], but in [13] it was shown that the square wave can be more efficient than the Haar wavelets in function reconstruction.

Moreover, when using the square wave basis, the compression algorithm is twice the complexity of the FFT, but the synthesis (decompression) is very simple [13]. Finally, as our approach has no limitation on the type of bases, one could get even better results employing different bases at different frequencies. For example using the sinusoid, the triangular wave or the sawtooth for the lower frequencies and square waves for the higher components.

We can also exploit the nonlinear properties of nonorthogonal bases to reveal the characteristics of complex systems. For example Vibrational Analysis is a technique for the study of the intrinsic dynamic characteristics of mechanical systems as a way to rate the state of their health. Even relatively simple mechanical systems have complex vibrational signatures, due to nonlinearities, that can be used to predict a potential failure. Up to now the tools in this field are mainly based on the FFT, while an analysis based on complex waveforms is better suited to the study of nonlinear phenomena and can result more effective in revealing the system behavior.

Other fields that could benefit from the increased degree of freedom in the basis choice are the signal (or pattern) recognition and cryptography. Usually, for recognition, an FFT of the signal is executed and the result is compared to a reference spectrum. Having the possibility of using any basis, one can perform the analysis employing a reference signal as a basis.

Another interesting application is in the field of noise suppression or signal separation, or as a way of building matched filters. If one analyzes a signal using as a basis one of the expected components, it is easy to separate noise from signal, or two signals that have been added, if they are at different repetition rates, essentially exploiting their S-orthogonality. As an example let us consider the signals in Fig. 6; at the top there is a square wave with an arbitrary phase, at the bottom there is the same square wave with some higher frequency noise added. When we analyze this resulting signal using the square wave as a basis, we obtain the generalized amplitude spectrum of Fig. 7, the phase is not plotted. Here, as we used the square wave as a basis, the first component of the spectrum completely represents the square wave at the fundamental frequency. The remaining part of the spectrum is relative to the decomposition of the noise in terms of series of square waves. So, from the first component one can reconstruct the original square wave (which is the equivalent of a "generalized" low pass filtering, in the square wave frequency domain). Whereas, applying a high pass filtering, i.e. keeping all the components of the spectrum except the first one, we can separate the noise portion as shown in Fig. 8. Moreover we could change the basis at any frequency to further exploit the technique. Note that this operation would be impossible using the standard sinusoidal filters, as the square wave has a (sinusoidal) spectrum that extends to infinite. This same technique can be used to surmount the limitations of the methods based on High Order Spectrum [11] (bispectrum and trispectrum) for frequency estimation.

## 6. Conclusions

The paper is intended as a first map of a huge new territory waiting to be explored. The objective was to show that a new decomposition method can be devised that works both in vectorial and polar coordinates. The theory generalizes the Fourier theorem to the case of nonorthogonal bases. Consequently a new paradigm seems opportune, to grant equal dignity to all $L^2$ bases and, even more, a new concept of "*generalized frequency domain*" should be introduced that is not limited to sinusoidal waveforms. Until now the sinusoidal (or complex exponential) basis has been considered as the center of the function space. Now we must recognize that that there is no center in that space. Any basis is equivalent, and we can pass from one representation to the other with an algorithm, essentially a change of coordinates. And, as it is done in other fields, one should always choose the coordinate frame that is more effective for the solution of the given problem.

The new tool has a number of advantages over the usual ones. Almost any function or collection of functions can be selected as a basis. This decomposition preserves some of the properties of the Fourier theorem so extending its applicability. With this approach the basis coincides with the *eigenfunction* of the corresponding system, this in turn paves the way for the development of new mathematical models of nonlinear phenomena. The algorithm is intuitive and very efficient.

The new tool, apart from simplifying the computation, can lead to interesting developments and applications, from the purely mathematical ones to immediate technical advances, and some example applications were given.

## *Acknowledgements*

The author wishes to thank the Italian Ministry of Foreign Affairs for having indirectly contributed to this work.